\documentclass[aop,MSNbibl,seceqn,secthm,dvips]{arximspdf}

\doi{10.1214/12-AOP789} %kopijuoti is PTS
\volume{41}
\issue{2}
\pubyear{2013}
\firstpage{1030}
\lastpage{1054}

\makeatletter
\newcommand{\rrvert}{\vert}
\newcommand{\llvert}{\vert}
\newproclaim{definition}[thm]{Definition}
\newtheorem{lemma}[thm]{Lemma}

\def\cal{\mathcal}
\def\al{{\alpha}}
\def\la{{\lambda}}
\def\si{{\sigma}}
\def\De{{\Delta}}
\def\Ga{{\Gamma}}
\def\Om{{\Omega}}

\def\CC{\mathbb C}
\def\EE{\mathbb E}
\def\HH{\mathbb H}
\def\NN{\mathbb N}
\def\RR{\mathbb R}
\def\ZZ{\mathbb Z}
\def\cB{{\cal B}}
\def\cF{{\cal F}}
\def\cG{{\cal G}}
\def\cS{{\cal S}}
\def\cU{{\cal U}}
\def\cX{{\cal X}}

\makeatother

\begin{document}
\begin{frontmatter}

\title{Super-Brownian motion as the unique strong solution to an SPDE}
\runtitle{Unique solution to SPDE}

\begin{aug}
\author[A]{\fnms{Jie} \snm{Xiong}\corref{}\thanksref{t2}\ead[label=e1]{jxiong@math.utk.edu}\ead[label=u1,url]{http://www.math.utk.edu/\textasciitilde jxiong}}
\thankstext{t2}{Supported in part by NSF DMS-09-06907.}
\runauthor{J. Xiong}
\affiliation{University of Tennessee and University of Macau}
\address[A]{Department of Mathematics\\
University of Tennessee\\
Knoxville, Tennessee 37996-1300\\
USA\\
\printead{e1}\\
\printead{u1}} %adresu isvedimo komanda gale!
\end{aug}

% HISTORY:
\received{\smonth{9} \syear{2010}}
\revised{\smonth{3} \syear{2012}}

% ABSTRACT
%
\begin{abstract}
A stochastic partial differential equation (SPDE)
is derived for super-Brownian motion regarded as a distribution
function valued process. The strong uniqueness for the solution to
this SPDE is obtained by an extended Yamada--Watanabe argument.
Similar results are also proved for the Fleming--Viot process.
\end{abstract}

% KEYWORDS


\end{frontmatter}

%s1 #&#
\section{Introduction}

Let $(\Om,\cF,P,\cF_t)$ be a stochastic basis satisfying the usual
conditions. Namely, $(\Om,\cF,P)$ is a probability space, and
$\{\cF_t\}$ is a family of nondecreasing right-continuous
sub-$\si$-fields of $\cF$ such that $\cF_0$ contains all $P$-null
subsets of $\Om$. Let $W$ be an $\cF_t$-adapted space--time white
noise random measure on $\RR_+\times U$ with intensity measure
$ds\,\la(da)$, where $(U,\cU,\la)$ is a measure space. We consider the
following stochastic partial differential equation (SPDE): for
$t\in\RR_+$ and $y\in\RR$,
%
%e1.1 #&#
\begin{equation}
\label{eq0623a} {u}_t(y)=F(y)+\int^t_0
\int_U G\bigl(a,y,{u}_s(y)\bigr){W}(ds \,da)+\int
^t_0\frac12\De{u}_s(y)\,ds,
\end{equation}
where $F$ is a
real-valued measurable function on $\RR$, $G\dvtx U\times\RR^2\to\RR$
satisfies the following conditions: there is a constant $K>0$ such
that for any $u_1, u_2, u,  y\in\RR$,
%
%e1.2 #&#
\begin{equation}
\label{eq0831d} \int_U\bigl|G(a,y,u_1)-G(a,y,u_2)\bigr|^2
\la(da)\le K|u_1-u_2|
\end{equation}
and
%
%e1.3 #&#
\begin{equation}
\label{eq0922a} \int_U\bigl|G(a,y,u)\bigr|^2\la(da)\le K
\bigl(1+|u|^2\bigr).
\end{equation}

We first give the definition for the solution to SPDE
(\ref{eq0623a}). To this end, we need to introduce the following
notation. For $i\in\NN\cup\{0\}$, let $\cX_i$ be the Hilbert space
consisting of all functions $f$ such that $f^{(k)}\in
L^2(\RR,e^{-|x|}\,dx)$, where $f^{(k)}$ denotes the $k$th order
derivative in the sense of generalized functions. We refer the
reader to Section 2.1 of Chapter 1 in the book of Gel'fand and
Shilov~\cite{GS} for a precise definition of such derivatives. We
shall denote $f^{(0)}=f$. The Hilbert norm $\|f\|_i$ is defined as
\[
\|f\|^2_i\equiv\sum_{k=0}^i
\int_\RR f^{(k)}(x)^2e^{-|x|}\,dx<
\infty.
\]
We denote the corresponding inner product by $\langle\cdot,\cdot
\rangle_i$.
Let $C^\infty_0(\RR)$ be the collection of functions which has
compact support and derivatives of all orders.

%de1.1 #&#
\begin{definition}
Suppose that $F\in\cX_0$. A continuous $\cX_0$-valued process
$\{{u}_t\}$ on a stochastic basis is a weak solution to SPDE
(\ref{eq0623a}) if there exists a space--time white noise ${W}$ such
that for any $t\ge0$ and $f\in C^\infty_0(\RR)$, we have
%
%e1.4 #&#
\begin{eqnarray}
\label{eq0820a} \langle{u}_t,f\rangle&=&\langle F,f\rangle+\int
^t_0\biggl\langle {u}_s,\frac12\De
f\biggr\rangle \,ds
\nonumber
\\[-8pt]
\\[-8pt]
\nonumber
&&{}+\int^t_0\int_\RR\int
_U G\bigl(a,y,{u}_s(y)\bigr)f(y)\,dy{W}(ds \,da)\qquad
\mbox{a.s.}
\end{eqnarray}
Here let
$\langle f,g\rangle=\int_\RR f(x)g(x)\,dx$ whenever the integral is
well-defined.

SPDE (\ref{eq0623a}) has a strong solution if for any space--time
white noise $W$ on stochastic basis $(\Om,\cF,P,\cF_t)$, there
exists a continuous $\cX_0$-valued $\cF_t$-adapted process
$\{{u}_t\}$ such that (\ref{eq0820a}) holds for all $f\in
C^\infty_0(\RR)$.
\end{definition}

The first main result of this article is presented as follows.
%
%th1.2 #&#
\begin{thm}\label{thm0618a}
Suppose that conditions (\ref{eq0831d}) and (\ref{eq0922a})
hold. If $F\in\cX_0$, then SPDE (\ref{eq0623a}) has a strong
solution $(u_t)$ satisfying
%
%e1.5 #&#
\begin{equation}
\label{eq0818a} \EE\sup_{0\le t\le T}\|u_t\|^2_0<
\infty,
\end{equation}
and any two
solutions satisfying this condition will coincide.
\end{thm}

The idea for the proof of the uniqueness part of Theorem
\ref{thm0618a} is outlined as follows. When the solution to SPDE
(\ref{eq0623a}) is $\cX_1$-valued, that is, $u_t(x)$ is differentiable
in $x$, we establish its connection to a backward doubly stochastic
differential equation (BDSDE). When the driving noise is finite
dimensional, the coefficients are Lipschitz, and the solution of the
SPDE is differentiable in $x$ up to order 2, this connection was
established by Pardoux and Peng~\cite{PP}. We will use a smoothing
approximation to achieve such a connection for the current
non-Lipschitz setting. The Yamada--Watanabe (cf.~\cite{YW}) argument
to the BDSDE is applied to establish the uniqueness of the solution.
As a consequence, SPDE (\ref{eq0623a}) has at most one solution in
the class of spatially differentiable solutions. In fact, the
uniqueness in this smaller space is sufficient for applications to
super-Brownian motions and Fleming--Viot processes.

The goal of Theorem~\ref{thm0618a} is to prove uniqueness in the set
of $\cX_0$-valued processes. The proof of this case is inspired by
that of the $\cX_1$-valued process. It uses a detailed estimate of
the spatial derivative term in the equation satisfied by the
smoothing approximation of the solutions.

The main motivation of the above result is its applications to many
measure-valued processes, from which three are stated here. At the
end of this section, other possible applications will be outlined,
while their presentations will appear in forthcoming publications.

Super-Brownian motion (SBM), also called the Dawson--Watanabe
process, has been studied by many authors since the pioneering work
of Dawson~\cite{D} and Watanabe~\cite{W}. It is a measure-valued
process arising as the limit for the empirical measure process of a
branching particle system. It has been proved that this process
satisfies a martingale problem (MP), whose uniqueness is
established by the nonlinear partial differential equation satisfied
by its log-Laplace transform. Denote SBM by $(\mu_t)$. When the
state space is $\RR$, for each $t$ and almost all ${\omega}$, the measure
$\mu_t$ has density with respect to the Lebesgue measure, and this
density-valued process $v_t$ satisfies the following nonlinear
SPDE:
%
%e1.6 #&#
\begin{equation}
\label{eq2a5} \partial_tv_t(x)=\tfrac12\De
v_t(x)+\sqrt{v_t(x)}\dot{B}_{tx},
\end{equation}
where $B$ is the space--time white noise on $\RR_+\times\RR$. This
SPDE was derived and studied independently by Konno and Shiga
\cite{KS} and Reimers~\cite{Rei}. The uniqueness of the solution to
SPDE (\ref{eq2a5}) is only proved in the weak sense using that of
the MP.

Many attempts have been made toward proving the strong uniqueness
for the solution to (\ref{eq2a5}). The main difficulty is the
non-Lipschitz coefficient in front of the noise. Some progress has
been made by relaxing the form of the SPDE. When the space $\RR$ is
replaced by a single point, (\ref{eq2a5}) becomes an SDE which is
the Feller's diffusion $dv_t=\sqrt{v_t}\,dB_t$ whose uniqueness is
established using the Yamada--Watanabe argument. When the random
field $B$ is colored in space and white in time, the strong
uniqueness of the solution to the SPDE (\ref{eq2a5}) with
$\sqrt{v_t(x)}$ replaced by a function of $v_t(x)$ was obtained by
Mytnik et al.~\cite{MPS} under suitable conditions. When $B$ is a
space--time white noise, Mytnik and Perkins~\cite{MP} prove pathwise
uniqueness for multiplicative noises of the form
$\si(x,v_t(x))\dot{B}_{tx}$, where $\si$ is H\"older continuous of
index $\al>\frac34$ in the solution variable. In particular, their
results imply that the SPDE
%
%e1.7 #&#
\begin{equation}
\label{eq0820b} \partial_tv_t(x)=\tfrac12\De
v_t(x)+\bigl|v_t(x)\bigr|^\al\dot{B}_{tx}
\end{equation}
has a pathwise unique solution when $\al>
\frac34$. Some negative results have also been achieved. When signed\vspace*{1pt}
solutions are allowed, Mueller et al.~\cite{MMP} give a
nonuniqueness result when $\frac12\le\al<\frac34$. For SPDE
(\ref{eq0820b}) restricted to nonnegative solutions, Burdzy
et al.~\cite{BMP} show a nonuniqueness result for
$0<\al<\frac12$.

In this paper, we approach this problem from a different point of
view. Instead of considering the equation for the density-valued
process, we study the SPDE satisfied by the ``distribution''
function-valued process. That is, we define the ``distribution''
function-valued process ${u}_t$,
%
%e1.8 #&#
\begin{equation}
\label{eq0831a} {u}_t(y)=\int^y_0
\mu_t(dx)\qquad \forall y\in\RR.
\end{equation}
Notice that ${u}_t(y)$ is differentiable
in $y$. Here $u_t$ is referred to as the corresponding distribution
function of $\mu_t$, although $\mu_t$ is not necessarily a
probability measure. In addition, we take the integral starting from
0 instead of $-\infty$ to include the case of $\mu_t$ being an
infinite measure.

Inspired by Dawson and Li~\cite{DL}, we consider the following SPDE:
%
%e1.9 #&#
\begin{equation}
\label{eq0609a} {u}_t(y)=F(y)+\int^t_0
\int^{{u}_s(y)}_0{W}(ds \,da)+\int^t_0
\frac 12\De {u}_s(y)\,ds,
\end{equation}
where $F(y)=\int^y_0\mu_0(dx)$ is the
distribution function of $\mu_0$, $W$ is a white noise random
measure on $\RR_+\times\RR$ with intensity measure $ds \,da$. The
authors of~\cite{DL} considered equation (\ref{eq0609a}) with
$\frac12\De$ replaced by the bounded operator $A$ given by
\[
Af(x)=\bigl({\gamma}(x)-f(x)\bigr)b,
\]
where $b$
is a constant and ${\gamma}$ is a fixed function. We prove that the
solution of~(\ref{eq0609a}) is indeed the distribution
function-valued process corresponding to an SBM. The strong
uniqueness for the solution to (\ref{eq0609a}) is then obtained by
applying Theorem~\ref{thm0618a} to the current setup. This result
provides a new proof of uniqueness in law for SBM.

%th1.3 #&#
\begin{thm}\label{thm4sbm}
Let $\{\mu_t\}$ be an SBM and $F\in\cX_0$, where
$F(y)\equiv\int^y_0\mu_0(dx)$, $\forall y\in\RR$. If $\{u_t\}$ is
the corresponding distribution function defined by (\ref{eq0831a}),
then it is possible to define a white noise $W$ on an extension of
the stochastic basis so that $\{u_t\}$ is the unique solution to the
SPDE (\ref{eq0609a}).

On the other hand, if $\{u_t\}$ is a weak solution to the SPDE
(\ref{eq0609a}) with $F\in\cX_0$ being nondecreasing, then
$\{\mu_t\}$ is an SBM.
\end{thm}

The definition of the extension of a stochastic basis and random
variables on the basis can be found in the book of Ikeda and
Watanabe~\cite{IW}. We refer the reader to Definition 7.1 on page 89
in~\cite{IW} for details. Here we only remark that the original SBM
remains an SBM on the extended stochastic basis.

Because of the difference in driving noise, the uniqueness of the
solution to SPDE (\ref{eq0609a}) does not imply that of SPDE
(\ref{eq2a5}). In fact, the noise $W$ in (\ref{eq0609a}) is
constructed using the noise $B$ and the solution $v_t$ in
(\ref{eq2a5}). We also note that our uniqueness of the solution to
SPDE (\ref{eq0609a}) does not contradict the nonuniqueness result
of~\cite{MMP} for the case of $\al=\frac12$, since signed solutions
are allowed in~\cite{MMP}. Let $v_t(x)$ be a (signed) solution to
(\ref{eq0820b}) with $\al=\frac12$. Then
\[
u_t(x)=\int^x_0v_t(y)\,dy
\]
does not satisfy [unless $v_t(x)$ is nonnegative] SPDE
(\ref{eq0609a}) because the quad\-ratic variation of the martingale
\[
\int^t_0\int^x_0\bigl|v_s(y)\bigr|^{1/2}B(ds
\,dy)
\]
is
\[
\int^t_0\int^x_0\bigl|v_s(y)\bigr|\,dy
\,ds\neq\int^t_0\bigl|u_s(x)\bigr|\,ds.
\]

Similarly, we consider another very important measure-valued
process: the Fleming--Viot (FV) process. We demonstrate that the
SPDE
%
%e1.10 #&#
\begin{equation}
\label{eq0609b}\qquad {u}_t(y)=F(y)+\int^t_0
\int^1_0\bigl(1_{a\le{u}_s(y)}-{u}_s(y)
\bigr){W}(ds \,da)+\int^t_0\frac12
\De{u}_s(y)\,ds
\end{equation}
can be used to
characterize the distribution function-valued process determined by
the FV process, where $W$ is a white noise random measure on
$\RR_+\times[0,1]$, with intensity measure $ds \,da$. Uniqueness of
the solution to SPDE (\ref{eq0609b}) is the second application of
Theorem~\ref{thm0618a}. Observe that this result provides a new
proof of uniqueness in law for FV process.

%th1.4 #&#
\begin{thm}\label{thm4fv}
Let $\{\mu_t\}$ be an FV process and
\[
u_t(y)=\mu_t\bigl((-\infty,y]\bigr)\qquad \forall y\in\RR.
\]
Let $F=u_0\in\cX_0$. Then it is possible to define a white noise
$W$ on an extension of the stochastic basis so that $\{u_t\}$ is the
unique solution to SPDE (\ref{eq0609b}).

On the other hand, if $\{u_t\}$ is a solution to SPDE
(\ref{eq0609b}) with $F\in\cX_0$ being the distribution of a
probability measure on $\RR$, then $\{\mu_t\}$ is an FV process.
\end{thm}

The third application of Theorem~\ref{thm0618a} is for the SPDE
driven by colored noise. More precisely, we consider the following
SPDE:
%
%e1.11 #&#
\begin{equation}
\label{eq0901a} du_t(x)=\tfrac12\De u_t(x)\,dt+
\sqrt{u_t(x)}B(x,dt),
\end{equation}
where $B$ is a Gaussian noise on $\RR\times\RR_+$ with covariance
function $\phi$ in space, that is,
\[
\EE B(x,dt)B(y,dt)=\phi(x,y)\,dt\qquad \forall x,y\in\RR.
\]

%th1.5 #&#
\begin{thm}\label{thm0904a}
Suppose $u_0\in\cX_0$ is fixed, and $\phi$ is bounded. Then SPDE (\ref{eq0901a}) has at most one solution.
\end{thm}

Such a result was obtained by Viot~\cite{V} when the state space is
bounded. The unbounded state space case was shown in~\cite{MPS}. We
reprove the result of~\cite{MPS} as an application of Theorem
\ref{thm0618a}. Mytnik, Perkins and Sturm~\cite{MPS} also consider
the case of singular
covariance; however, Theorem~\ref{thm0618a} does not apply to this
case.

The rest of this paper is organized as follows. In Section
\ref{sec2}, we establish the existence of a solution to SPDE
(\ref{eq0623a}). Section~\ref{sec3} introduces the BDSDE and gives a
Yamada--Watanabe type criteria for such equation. It also illustrates
the connection between the SPDE and the BDSDE. As a consequence,
uniqueness for the solution of the SPDE when the solutions are
restricted to those with first order partial derivative in the
spatial variable. We refine in Section~\ref{sec4} the uniqueness
proof of Section~\ref{sec3} without the spatial differentiability
condition. Finally, Section~\ref{sec5} applies the uniqueness result
for SPDE (\ref{eq0623a}) to three important measure-valued
processes.

We use $\mu(f)$ or $\langle\mu,f\rangle$ to denote the integral of
a function
$f$ with respect to the measure $\mu$. The letter $K$ stands for a
constant whose value can be changed from place to place.
$\partial_x$ is used to denote the partial derivative with respect
to the variable $x$ if the notation $\nabla$
is ambiguous.

We conclude this section by mentioning other possible applications
of the idea developed in this article. The first is to consider
measure-valued processes with interaction among individuals in the
system. This interaction may come from the drift and diffusion
coefficients which govern the motion of the individuals. It may also
come from the branching and immigration mechanisms. This extension
will appear in a joint work of Mytnik and Xiong~\cite{MX}. The
second possible application is to consider other type of nonlinear
SPDEs, especially those where the noise term involves the spatial
derivative of the solution. This extension will appear in a joint
work of Gomez et al.~\cite{GLMWX}. Finally, studying
measure-valued processes by using SPDE methodology will have the
advantage of utilizing the rich collection of tools developed in the
area of SPDEs. For example, the large deviation principle (LDP) for
some measure-valued processes, including FV process and the SBM, can
be established. As is well known, LDP for general FV process is a
long standing open problem (some partial results were obtained by
Dawson and Feng~\cite{DF1,DF2}, and Feng and Xiong~\cite{FX}
for neutral FV processes, and Xiang and Zhang~\cite{XZ} for the case
when the mutation operator tends to 0). This application will be
presented in a joint work of Fatheddin and Xiong~\cite{FatX}.

It was pointed out to me by two referees and by Leonid Mytnik that
Theorem~\ref{thm0618a} can be proved using the Yamada--Watanabe
argument directly to the SPDE without introducing the BDSDE. One of
the advantages of the current backward framework is that the term
involving the Laplacian operator gets canceled when
It\^o--Pardoux--Peng formula is applied. Furthermore, as one of the
referees pointed out, ``it is quite possible that the BDSDE idea
will have something to offer in other natural interacting models.''
In fact, in~\cite{GLMWX}, the BDSDE idea is used to get the
uniqueness for the solution to an SPDE where the noise term involves
the spatial derivative of the solution. This term actually helped us
in the proof of the uniqueness of the solution. To the best of my
knowledge, the direct Yamada--Watanabe argument to such an equation
cannot be easily implemented in this case.

%s2 #&#
\section{Existence of solution to SPDE}\label{sec2}

In this section, we consider the existence of a solution to SPDE
(\ref{eq0623a}).

Note that the definition of weak solution to (\ref{eq0623a}) is
equivalent to the following mild formulation:
%
%e2.1 #&#
\begin{equation}
\label{eq0615a}\quad u_t(y)=T_tF(y)+\int^t_0
\int_U\int_\RR p_{t-s}(y-z)G
\bigl(a,z,u_s(z)\bigr)\,dz W(ds \,da),
\end{equation}
where $T_t$ is the Brownian semigroup, which is for any $f\in\cX_0$,
\[
T_t f(x)=\int_\RR p_t(x-y)f(y)\,dy\quad
\mbox{and}\quad p_t(x)=\frac{1}{\sqrt{2\pi t}}\exp\biggl(-\frac{x^2}{2t}
\biggr).
\]

Before constructing a solution to (\ref{eq0615a}), we prove the
semigroup property for the family $\{T_t\}$ to be used in later
sections.

%le2.1 #&#
\begin{lemma}\label{lem0225a}
$\{T_t\dvtx t\ge0\}$ is a strongly continuous semigroup on $\cX_0$.
\end{lemma}
\begin{pf} Let $K_t$ be the function given by
\[
K_t^2=\int_\RR
e^{t|z|}p_1(z)\,dz<\infty\qquad \forall t\ge0.
\]
It is easy to show that for any $f\in\cX_0$, we have
%
%e2.2 #&#
\begin{equation}
\label{eq0225a} \|T_tf\|_0\le K_t\|f
\|_0.
\end{equation}
Thus, $\{T_t, t\ge0\}$ is a family of bounded linear operators on
$\cX_0$. The semigroup property is not difficult to verify. We now
focus on this semigroup's strong continuity.

For any $f\in C_b(\RR)\cap\cX_0$, it follows from the dominated
convergence theorem that as $t\to0$,
\[
\|T_tf-f\|^2_0\le\int_\RR
\biggl\llvert \int_\RR \bigl(f(x+tz)-f(x)
\bigr)p_1(z)\,dz\biggr\rrvert^2 e^{-|x|}\,dx\to0.
\]
In general, for $f\in\cX_0$, we take a sequence
$f_n\in C_b(\RR)\cap\cX_0$ such that $\|f_n-f\|_0\to0$ as
$n\to\infty$. Then
\[
\|T_tf-f\|_0\le K_t\|f_n-f
\|_0+\|T_tf_n-f_n
\|_0,
\]
which implies $T_tf\to f$ in $\cX_0$ as $t\to0$.
\end{pf}

In addition, we define operators $T^U_t$ on the Hilbert space
$\cX_0\otimes L^2(U,\la)=L^2(\RR\times U,e^{-|x|}\,dx\,\la(da))$ as
\[
T^U_tg(a,x)=\int_\RR
p_t(x-y)g(a,y)\,dy\qquad \forall t\ge0.
\]
By the same argument as in the proof of Lemma~\ref{lem0225a}, we
have the following result.

%le2.2 #&#
\begin{lemma}\label{lem0225b}
$\{T^U_t\dvtx t\ge0\}$ is a strongly continuous semigroup on
$\cX_0\otimes L^2(U,\la)$. Furthermore, for any $g\in\cX_0\otimes
L^2(U,\la)$,
%
%e2.3 #&#
\begin{equation}
\label{eq0225b} \bigl\|T^U_tg\bigr\|_{\cX_0\otimes L^2(U,\la)}\le
K_t\|g\|_{\cX_0\otimes
L^2(U,\la)}.
\end{equation}
\end{lemma}

Now, we come back to the construction of a solution to
(\ref{eq0615a}). Define a sequence of approximations by:
$u^0_t(y)=F(y)$ and, for $n\ge0$,
%
%e2.4 #&#
\begin{equation}
\label{eq0615b} \qquad u^{n+1}_t(y)=T_tF(y)+\int
^t_0\int_U\int
_\RR p_{t-s}(y-z)G\bigl(a,z,u^n_s(z)
\bigr)\,dz W(ds \,da).
\end{equation}

Let
\[
J(x)=\int_\RR e^{-|y|}\rho(x-y)\,dy,
\]
where $\rho$ is the mollifier given by
\[
\rho(x)=K\exp\bigl(-1/\bigl(1-x^2\bigr)\bigr)1_{|x|<1},
\]
and $K$ is a constant such that $\int_\RR\rho(x)\,dx=1$. Then, for any
$m\in\ZZ_+$, there are constants $c_m$ and $C_m$ such that
\[
c_me^{-|x|}\le J^{(m)}(x)\le C_me^{-|x|}\qquad
\forall x\in\RR;
\]
cf. Mitoma~\cite{Mit}, (2.1). We may and will replace $e^{-|x|}$
by $J(x)$ in the definition of space $\cX_i$.

%le2.3 #&#
\begin{lemma}
For any $p\ge1$ and $T>0$, there exists a constant $K_1=K_1(p,T)$
such that for any $n\ge0$,
%
%e2.5 #&#
\begin{equation}
\label{eq0818b} \EE\sup_{t\le T}\bigl\|u^n_t
\bigr\|^{2p}_0\le K_1.
\end{equation}
\end{lemma}
\begin{pf} We proceed by adapting the idea of Kurtz and Xiong~\cite{KX}.
Smoothing out if necessary, we may and will assume that
$u^{n+1}_t\in\cX_2$. By It\^o's formula, it is easy to show that,
for any $f\in C^\infty_0(\RR)$,
%
%e2.6 #&#
\begin{eqnarray}
\label{eq0615c}\qquad \bigl\langle{u}^{n+1}_t,f\bigr
\rangle_0&=&\langle F,f\rangle_0+\int^t_0
\biggl\langle\frac12\De{u}^{n+1}_s, f\biggr
\rangle_0\,ds
\nonumber
\\[-8pt]
\\[-8pt]
\nonumber
&&{}+\int^t_0\int_\RR\int
_U G\bigl(a,y,{u}^n_s(y)
\bigr)f(y)J(y)\,dy{W}(ds \,da)\qquad \mbox{a.s.}
\end{eqnarray}

Applying It\^o's formula to (\ref{eq0615c}) gives
\begin{eqnarray*}
&&\bigl\langle{u}^{n+1}_t,f\bigr\rangle^2_0
\\
&&\qquad=\langle F,f\rangle_0^2+\int^t_0
\bigl\langle{u}^{n+1}_s,f\bigr\rangle_0\bigl
\langle \De{u}^{n+1}_s,f\bigr\rangle_0\,ds
\\
&&\qquad\quad{}+\int^t_0\int_U\biggl(
\int_\RR G\bigl(a,y,{u}^n_s(y)
\bigr)f(y)J(y)\,dy\biggr)^2\la(da) \,ds
\\
&&\qquad\quad{}+\int^t_0\int_U 2\bigl
\langle{u}^{n+1}_s,f\bigr\rangle_0\int
_\RR G\bigl(a,y,{u}^n_s(y)
\bigr)f(y)J(y)\,dy{W}(ds \,da).
\end{eqnarray*}
Summing on $f$
over a complete orthonormal system (CONS) of $\cX_0$, we have
\begin{eqnarray*}
\bigl\|{u}^{n+1}_t\bigr\|_0^2&=&\|F
\|_0^2+\int^t_0\bigl
\langle{u}^{n+1}_s,\De{u}^{n+1}_s
\bigr\rangle_0\,ds
\\
&&{}+\int^t_0\int_U\int
_\RR G\bigl(a,y,{u}^n_s(y)
\bigr)^2J(y)\,dy \,\la(da) \,ds
\\
&&{}+\int^t_0\int_U 2\bigl
\langle{u}^{n+1}_s, G\bigl(a,\cdot,{u}^n_s(
\cdot)\bigr)\bigr\rangle_0{W}(ds \,da).
\end{eqnarray*}
It\^o's
formula is again applied to obtain
%
%e2.7 #&#
\begin{eqnarray}
\label{eq0615d}
&&\bigl\|{u}^{n+1}_t\bigr\|_0^{2p}
\nonumber\\
&&\qquad=\|F\|_0^{2p}+\int^t_0p
\bigl\|u^{n+1}_s\bigr\|^{2(p-1)}_0\bigl\langle
{u}^{n+1}_s,\De{u}^{n+1}_s\bigr
\rangle_0\,ds
\nonumber
\\
&&\qquad\quad{}+\int^t_0p\bigl\|u^{n+1}_s
\bigr\|^{2(p-1)}_0\int_U\int
_\RR G\bigl(a,y,{u}^n_s(y)
\bigr)^2J(y)\,dy\, \la(da) \,ds
\\
&&\qquad\quad{}+\int^t_0p\bigl\|u^{n+1}_s
\bigr\|^{2(p-1)}_0\int_U 2\bigl
\langle{u}^{n+1}_s, G\bigl(a,\cdot,{u}^n_s(
\cdot)\bigr)\bigr\rangle_0{W}(ds \,da)
\nonumber
\\
&&\qquad\quad{}+2p(p-1)\int^t_0\bigl\|u^{n+1}_s
\bigr\|^{2(p-2)}_0\int_U \bigl
\langle{u}^{n+1}_s, G\bigl(a,\cdot,{u}^n_s(
\cdot)\bigr)\bigr\rangle_0^2\la(da) \,ds.
\nonumber
\end{eqnarray}
Note that,
for $u\in\cX_1$,
\[
\int_\RR u(x)u'(x)J'(x)\,dx =-
\int_\RR u(x) \bigl(u'(x)J'(x)+u(x)J''(x)
\bigr)\,dx,
\]
which implies that
\begin{eqnarray*}
-\int_\RR u(x)u'(x)J'(x)\,dx&=&
\frac12\int_\RR u(x)^2J''(x)\,dx
\le K_2\int_\RR u(x)^2J(x)\,dx
 =K_2\|u\|_0^2.
\end{eqnarray*}
Therefore,
\begin{eqnarray*}
\langle u,\De u\rangle_0&=&\int_\RR
u''(x)u(x)J(x)\,dx
\\
&=&-\int_\RR u'(x) \bigl(u'(x)J(x)+u(x)J'(x)
\bigr)\,dx
\\
&\le&K_2\|u\|_0^2.
\end{eqnarray*}
By using the Burkholder--Davis--Gundy inequality on (\ref{eq0615d}),
\begin{eqnarray*}
- \EE\sup_{s\le t}\bigl\|{u}^{n+1}_s\bigr\|_0^{2p}&
\le&\|F\|_0^{2p}+pK_2\int^t_0
\EE\bigl\|u^{n+1}_s\bigr\|^{2p}_0\,ds
\\
&&+K_3\int^t_0\EE\bigl(
\bigl\|u^{n+1}_s\bigr\|^{2(p-1)}_0\bigl(1+
\bigl\|{u}^n_s\bigr\|^2_0\bigr)\bigr)\,ds
\\
&&{}+K_4\EE\biggl(\int^t_0
\bigl\|u^{n+1}_s\bigr\|^{4p-2}_0\bigl(1+
\bigl\|{u}^n_s\bigr\|^2_0\bigr)\,ds
\biggr)^{1/2}.
\end{eqnarray*}
Hence,
\begin{eqnarray*}
f_{n+1}(t)&\equiv&\EE\sup_{s\le t}\bigl\|{u}^{n+1}_s
\bigr\|_0^{2p}
\\
&\le&\|F\|_0^{2p}+K_5\int
^t_0f_{n+1}(s)\,ds+K_6\int
^t_0f_n(s)\,ds\\
&&{}+\frac
{1}{2}f_{n+1}(t).
\end{eqnarray*}
Gronwall's inequality and an induction argument finish the
proof.
\end{pf}

We proceed to proving the tightness of $\{u^n\}$ in
$C([0,T]\times\RR)$. Denote
\[
v^n_t(y)=\int^t_0\int
_U\int_\RR p_{t-s}(y-z)G
\bigl(a,z,u^n_s(z)\bigr)\,dz W(ds \,da).
\]

%le2.4 #&#
\begin{lemma}
For any $p\ge1>\al$, there is a constant $K_1$ such that
%
%e2.8 #&#
\begin{equation}
\label{eq0825a} \EE\bigl|v^n_t(y_1)-v^n_t(y_2)\bigr|^{2p}
\le K_1e^{p(|y_1|\vee|y_2|)}|y_1-y_2|^{p\al}.
\end{equation}
\end{lemma}

\begin{pf} Denote the left-hand side of (\ref{eq0825a}) by $I$. It
follows from Burkhol\-der's inequality that there exists a constant
$K_2>0$ such that $I$ is bounded by
\[
K_2 \EE\biggl(\int^t_0\int
_U\biggl(\int_\RR
\bigl(p_{s}(y_1-z)-p_{s}(y_2-z)
\bigr)G\bigl(a,z,u^n_{t-s}(z)\bigr)\,dz\biggr)^2
\la(da) \,ds\biggr)^p.
\]
By H\"older's inequality,
\begin{eqnarray*}
I&\le&K_2\EE \biggl(\int^t_0\int
_U\int_\RR \bigl(p_{s}(y_1-z)-p_{s}(y_2-z)
\bigr)^2e^{|z|}\,dz
\\
&&\hspace*{66pt}{} \times\int_\RR G\bigl(a,z,u^n_{t-s}(z)
\bigr)^2e^{-|z|}\,dz\, \la(da) \,ds \biggr)^p.
\end{eqnarray*}
The linear growth condition (\ref{eq0922a}) and the estimate
(\ref{eq0818b}) is then applied to get
\begin{eqnarray*}
I&\le&K_2\EE \biggl(\int^t_0\int
_\RR \bigl(p_{s}(y_1-z)-p_{s}(y_2-z)
\bigr)^2e^{|z|}\,dz\\
&&\hspace*{55pt}{}\times
\int_\RR K\bigl(1+\bigl|u^n_{t-s}(z)\bigr|^2
\bigr)e^{-|z|}\,dz \,ds \biggr)^p
\\
&\le&K_3\biggl(\int^t_0\int
_\RR \bigl(p_{s}(y_1-z)-p_{s}(y_2-z)
\bigr)^2e^{|z|}\,dz\,ds\biggr)^p.
\end{eqnarray*}
Using the fact that
\[
\bigl|p_s(y_1)-p_s(y_2)\bigr|\le
Ks^{-1}|y_1-y_2|\qquad \forall s>0, y_1,
y_2\in\RR,
\]
we arrive at
\begin{eqnarray*}
I&\le&K_4\biggl(\int^t_0\int
_\RR s^{-\al}|y_1-y_2|^\al
\bigl(p_s(y_1-z)\vee p_s(y_2-z)
\bigr)^{2-\al}e^{|z|}\,dz\,ds\biggr)^p
\\
&\le&K_4\biggl(\int^t_0\int
_\RR s^{-\al}|y_1-y_2|^\al
p_s(z)^{2-\al}e^{|z|}\,dz\,ds\, e^{|y_1|\vee|y_2|}
\biggr)^p
\\
&\le&K_5\biggl(\int^t_0
s^{-\al} s^{-(1-\al)/2}\,ds\biggr)^pe^{p(|y_1|\vee|y_2|)}|y_1-y_2|^{p\al}
\\
&\le&K_1e^{p(|y_1|\vee|y_2|)}|y_1-y_2|^{p\al},
\end{eqnarray*}
which finishes the proof of (\ref{eq0825a}).
\end{pf}

Similarly, we can prove that
\[
\EE\bigl|v^n_{t_1}(y)-v^n_{t_2}(y)\bigr|^{2p}
\le K_1e^{p|y|/2}|t_1-t_2|^{p\al/2}.
\]

We are now ready to provide.

\begin{pf*}{Proof of Theorem~\ref{thm0618a} \textup{(Existence)}} By Kolmogorov's
criteria (cf. Corollary 16.9 in Kallenberg~\cite{Kallen}), for each
fixed $m$, the sequence of laws of $\{v^n_t(x)\dvtx (t,x)\in[0,T]\times[-m,m]\}$ on $\CC([0,T]\times[-m,m])$ is tight,
and hence, has a convergent subsequence. By the standard
diagonalization argument, there exists a subsequence
$\{v^{n_k}_t(x)\}$ which converges in law on
$\CC([0,T]\times[-m,m])$ for each $m$. Therefore, $\{v^{n_k}_t(x)\}$
converges in law on $\CC([0,T]\times\RR)$.

Let $v_t(x)$ be a limit point. For any $t_1<t_2$, it follows from
Fatou's lemma that
\begin{eqnarray*}
\EE\|v_{t_1}-v_{t_2}\|^{2p}_0&
\le&K_1\liminf_{k\to\infty}\EE\biggl(\int_\RR\bigl|v^{n_k}_{t_1}(x)-v^{n_k}_{t_2}(x)\bigr|^2e^{-|x|}\,dx
\biggr)^p
\\
&\le&K_2\liminf_{k\to\infty}\EE\int_\RR
\bigl|v^{n_k}_{t_1}(x)-v^{n_k}_{t_2}(x)\bigr|^{2p}e^{-(2/3)p|x|}\,dx
\\
&\le&K_3 \int_\RR e^{(1/2)p|x|}|t_1-t_2|^{p\al/2}e^{-(2/3)p|x|}\,dx
\\
&=&K_4|t_1-t_2|^{p\al/2}.
\end{eqnarray*}
By Kolmogorov's criteria again, we see that there is a version,
which we will take, such that $v_\cdot\in\CC([0,T],\cX_0)$ a.s.

Let $u_t(y)=T_tF(y)+v_t(y)$. Then, $u_\cdot\in\CC([0,T],\cX_0)$
a.s. The proof of $\{u_\cdot\}$ being a solution to SPDE
(\ref{eq0623a}) is standard. Here is a sketch and the reader is
referred to Sections 6.2 and 8.2 of Kallianpur and Xiong~\cite{KX1}
for two similar situations. First, by passing to the limit, we can
prove that for any $f\in C^\infty_0(\RR)$,
\[
M^f_t\equiv\langle{u}_t,f\rangle-\langle
F,f\rangle-\int^t_0\biggl\langle
{u}_s,\frac12\De f\biggr\rangle \,ds
\]
and
\begin{eqnarray*}
N^f_t&\equiv&\langle{u}_t,f
\rangle^2-\langle F,f\rangle^2-\int^t_0
\langle {u}_s,f\rangle \langle {u}_s,\De f\rangle \,ds
\\
&&{}-\int^t_0\int_U\biggl(
\int_\RR G\bigl(a,y,u_s(y)\bigr)f(y)\,dy
\biggr)^2\la(da)\,ds
\end{eqnarray*}
are martingales. It then follows that the quadratic
variation process of $M^f$ is given by
\[
\bigl\langle M^f\bigr\rangle_t=\int
^t_0\int_U\biggl(\int
_\RR G\bigl(a,y,u_s(y)\bigr)f(y)\,dy
\biggr)^2\la(da)\,ds.
\]
The martingale $M^f$ is then represented as
\[
M^f_t=\int^t_0\int
_\RR\int_U G\bigl(a,y,{u}_s(y)
\bigr)f(y)\,dy{W}(ds \,da)
\]
for a suitable random measure $W$ defined on a stochastic basis.
Consequently, $u_t$ is a weak solution to SPDE (\ref{eq0623a}).

Estimate (\ref{eq0818a}) follows
from (\ref{eq0818b}) and Fatou's lemma.
\end{pf*}

%s3 #&#
\section{Backward doubly SDE}\label{sec3}

This section is of interest on its own. It is inspirational for the
proof of the uniqueness part of Theorem~\ref{thm0618a}, which we
will present in the next section.

In this section, we study uniqueness of the solution to a BDSDE
whose coefficient is not Lipschitz, and the relationship between
this BDSDE and an SPDE whose coefficient is not Lipschitz. Because
of this non-Lipschitz property, the corresponding results of Pardoux
and Peng~\cite{PP} do not apply to the current BDSDE and SPDE. We
will adapt Yamada--Watanabe's argument to the present setup to obtain
uniqueness for the solution to the BDSDE and a smoothing
approximation to establish the connection between the BDSDE and the
SPDE. As an application, we obtain the uniqueness for the SPDE if
the solutions are restricted to those that are differentiable with
respect to the spatial variable.

Let $y\in\RR$ be fixed. We consider the following BDSDE with pair
$(Y_t,Z_t)$ as its solution:
%
%e3.1 #&#
\begin{equation}
\label{eq0831c}\quad Y_t=\xi+\int^T_t
\int_U G(a,y,Y_s)\tilde{W}(\hat{d}s \,da)-\int
^T_tZ_s\,dB_s,\qquad 0\le t
\le T,
\end{equation}
where $\xi$ is an $\cF^B_T$-measurable random variable, $G$
satisfies the H\"older continuity~(\ref{eq0831d}),
$\cF^B_T=\si(B_s\dvtx 0\le s\le T)$, $B$ is a Brownian motion and
$\tilde{W}$, independent of $B$, is a space--time white noise in
$\RR_+\times U$ with intensity measure $ds \,\la(da)$. The notation
$\tilde{W}(\hat{d}s \,da)$ stands for the backward It\^o integral (cf.
Xiong~\cite{X}), that is, in the Riemann sum approximating the
stochastic integral, we take the right endpoints instead of the
left ones.

%de3.1 #&#
\begin{definition}
The pair of processes $(Y_t,Z_t)$ is a solution to BDSDE
(\ref{eq0831c}) if they are $\cG_t$-adapted, $Y_\cdot\in
C([0,T],\RR)$ a.s., $\EE\int^T_0Z_s^2\,ds<\infty$ and for each
$t\in[0,T]$, identity (\ref{eq0831c}) holds a.s., where
$\cG_t=\si(\cF^B_t,\cG^1_t)$ and $\cG^1_t$ is a nonincreasing
family of $\si$-fields which is independent of $B$ and contains
\[
\cF^{\tilde{W}}_{t,T}=\si\bigl(\tilde{W}\bigl([r,T]\times A\bigr), r
\in[t,T], A\in\cB(\RR)\bigr).
\]
\end{definition}

Note that the family $\{\cG_t\}$ is not a filtration because it is
not increasing. We now state an It\^o type formula in the present
setting.

%le3.2 #&#
\begin{lemma}[(It\^o--Pardoux--Peng formula)]
Suppose that a process $y_t$ is given by
\[
y_t=\xi+\int^T_t\int
_U\al(s,a)\tilde{W}(\hat{d}s \,da)-\int^T_tz_s\,dB_s,
\]
where $\al\dvtx  [0,T]\times U\times\Om\to\RR$ is a $\cG_t$-adapted
random field, and
\[
\EE\int^T_0\int_U
\al(s,a)^2\la(da)\,ds+\EE\int^T_0z_s^2\,ds<
\infty.
\]
Then, for any $f\in C^2_b(\RR)$, we have
%
%e3.2 #&#
\begin{eqnarray}
\label{eq0824a} f(y_t)&=&f(\xi)+\int^T_t
\int_U f'(y_s) \al(s,a)\tilde{W}(
\hat{d}s \,da)-\int^T_tz_sf'(y_s)\,dB_s
\nonumber
\\[-8pt]
\\[-8pt]
\nonumber
&&{}+\frac12\int^T_t\int_U
f''(y_s) \al(s,a)^2 \,da \,ds-
\frac12 \int^T_tz^2_sf''(y_s)\,ds.
\end{eqnarray}
\end{lemma}
\begin{pf} Let $\{h_j\}$ be a CONS of $L^2(U,\cU,\la)$ and
\[
\tilde{W}^{h_j}_t=\int^t_0
\int_Uh_j(a)\tilde{W}(ds \,da), \qquad j=1,2,\ldots.
\]
Then, $\{\tilde{W}^{h_j}_t\}_{j=1,2,\ldots}$ are
independent Brownian motions. Let
\[
y^n_t=\xi+\sum^n_{j=1}
\int^T_t\bigl\langle\al(s,\cdot),h_j
\bigr\rangle_{L^2(U,\la)}\,\hat{d}\tilde{W}^{h_j}_s-\int
^T_tz_s\,dB_s,
\]
where $\langle\cdot,\cdot\rangle_{L^2(U,\la)}$ denotes the inner
product in
$L^2(U,\cU,\la)$, and $\hat{d}\tilde{W}^{h_j}_s$ means that the
stochastic integral is defined as backward It\^o integral.

Applying Lemma 1.3 of~\cite{PP} to $f(y^n_t)$ gives
\begin{eqnarray*}
f\bigl(y^n_t\bigr)&=&f(\xi)+\sum
^n_{j=1}\int^T_t
f'\bigl(y^n_s\bigr) \bigl\langle\al(s,
\cdot),h_j\bigr\rangle_{L^2(U,\la)}\,\hat{d}\tilde
{W}^{h_j}_s-\int^T_tz_sf'
\bigl(y^n_s\bigr)\,dB_s
\\
&&{}+\frac12\sum^n_{j=1}\int
^T_t f''
\bigl(y^n_s\bigr)\bigl\langle\al(s,\cdot),h_j
\bigr\rangle_{L^2(U,\la)}^2 \,ds-\frac12 \int^T_tz^2_sf''
\bigl(y^n_s\bigr)\,ds.
\end{eqnarray*}
Taking $n\to\infty$, we then finish the proof of the
It\^o--Pardoux--Peng formula (\ref{eq0824a}) under the present setup.
\end{pf}

Here is the main result of this section.

%th3.3 #&#
\begin{thm}\label{thm0901b}
Suppose that conditions (\ref{eq0831d}) and (\ref{eq0922a})
hold. Then\break BDSDE~(\ref{eq0831c}) has at most one solution.
\end{thm}
\begin{pf} Suppose that (\ref{eq0831c}) has two solutions
$(Y^i_t,Z^i_t),  i=1, 2$. Let $\{a_k\}$ be a decreasing positive
sequence defined recursively by
\[
a_0=1\quad\mbox{and}\quad \int^{a_{k-1}}_{a_k}z^{-1}\,dz=k,\qquad
k\ge1.
\]
Let $\psi_k$ be nonnegative continuous functions supported in
$(a_k,a_{k-1})$ satisfying
\[
\int^{a_{k-1}}_{a_k}\psi_k(z)\,dz=1\quad\mbox{and}\quad\psi_k(z)\le 2(kz)^{-1}\qquad \forall z\in\RR.
\]
Let
\[
\phi_k(z)=\int^{|z|}_0\,dy\int
^y_0\psi_k(x)\,dx\qquad \forall z\in \RR.
\]
Then, $\phi_k(z)\to|z|$ and $|z|\phi''_k(z)\le2k^{-1}$.

Since
%
%e3.3 #&#
\begin{eqnarray}
\label{eq0831e} Y^1_t-Y^2_t&=&
\int^T_t\int_U \bigl(G
\bigl(a,y,Y^1_s\bigr)-G\bigl(a,y,Y^2_s
\bigr)\bigr)\tilde{W}(\hat{d}s \,da)
\nonumber
\\[-9pt]
\\[-9pt]
\nonumber
&&{}-\int^T_t\bigl(Z^1_s-Z^2_s
\bigr)\,dB_s,
\end{eqnarray}
then by the
It\^o--Pardoux--Peng formula,
%
%e3.4 #&#
\begin{eqnarray}
\label{eq0902c}
&&\phi_k\bigl(Y^1_t-Y^2_t
\bigr)
\nonumber
\\[-2pt]
&&\qquad=\int^T_t\int_U
\phi'_k\bigl(Y^1_s-Y^2_s
\bigr) \bigl(G\bigl(a,y,Y^1_s\bigr)-G
\bigl(a,y,Y^2_s\bigr)\bigr)\tilde{W}(\hat{d}s \,da)
\nonumber
\\[-2pt]
&&\qquad\quad{}-\int^T_t\phi'_k
\bigl(Y^1_s-Y^2_s\bigr)
\bigl(Z^1_s-Z^2_s
\bigr)\,dB_s
\\[-2pt]
&&\qquad\quad{}+\frac12\int^T_t\int_U
\phi''_k\bigl(Y^1_s-Y^2_s
\bigr) \bigl(G\bigl(a,y,Y^1_s\bigr)-G
\bigl(a,y,Y^2_s\bigr)\bigr)^2\la(da) \,ds
\nonumber
\\[-2pt]
&&\qquad\quad{}-\frac12\int^T_t\phi''_k
\bigl(Y^1_s-Y^2_s\bigr)
\bigl(Z^1_s-Z^2_s
\bigr)^2\,ds.\nonumber
\end{eqnarray}
The sequence $\phi'_k$ being bounded and
$\EE\int^T_0|Z^1_s-Z^2_s|^2\,ds<\infty$ imply that the second term on
the right-hand side of (\ref{eq0902c}) is a square integrable
martingale, and hence, its expectation is 0. Moreover, by a parallel
argument, the expectation of the first term is also zero. Since the
last term is nonpositive, by taking expectation on both sides of
(\ref{eq0902c}), the following estimate is attained
\begin{eqnarray*}
&&\EE\phi_k\bigl(Y^1_t-Y^2_t
\bigr)
\\
&&\qquad\le\EE\frac12\int^T_t\int
_U\phi''_k
\bigl(Y^1_s-Y^2_s\bigr) \bigl(G
\bigl(a,y,Y^1_s\bigr)-G\bigl(a,y,Y^2_s
\bigr)\bigr)^2\la(da) \,ds
\\
&&\qquad\le K_1\EE\int^T_t
\phi''_k\bigl(Y^1_s-Y^2_s
\bigr)\bigl|Y^1_s-Y^2_s\bigr|\,ds
\\
&&\qquad\le K_2k^{-1}.
\end{eqnarray*}
Taking $k\to\infty$ and making use of Fatou's lemma, we have
\[
\EE\bigl|Y^1_t-Y^2_t\bigr|\le0.
\]
Therefore, $Y^1_t=Y^2_t$ a.s. Plugging back into (\ref{eq0831e}), we
can get
\[
\int^T_t\bigl(Z^1_s-Z^2_s
\bigr)\,dB_s=0 \qquad \mbox{a.s.}
\]
Hence,
$Z^1_t=Z^2_t$ a.s. for a.e. $t$, completing the proof.
\end{pf}

Finally, in this section, we establish a relationship between SPDEs
and BDSDEs under non-Lipschitz setup. To this end, we convert SPDE
(\ref{eq0623a}) to its backward version. For $T$ fixed, we define
the random field
\[
\tilde{u}_t(y)=u_{T-t}(y)\qquad\forall t\in[0,T], y\in\RR,
\]
and introduce the new noise $\tilde{W}$ by
\[
\tilde{W}\bigl([0,t]\times A\bigr)=W\bigl([T-t,T]\times A\bigr)\qquad \forall t\in
[0,T], A\in\cB(\RR).
\]
Then, $\tilde{u}_t$ satisfies backward SPDE given by
%
%e3.5 #&#
\begin{equation}
\label{eq2a6}\quad \tilde{u}_t(y)=F(y)+\int^T_t
\int_U G\bigl(a,y,\tilde{u}_s(y)\bigr)
\tilde{W}(\hat{d}s \,da)+\int^T_t\frac12\De
\tilde{u}_s(y)\,ds.
\end{equation}
It is
clear that SPDEs (\ref{eq0623a}) and (\ref{eq2a6}) have the same
uniqueness property. Specifically, if (\ref{eq0623a}) has a unique
strong solution, then so does (\ref{eq2a6}), and vice versa. Observe
that $\tilde{u}_t$ is $\cF^{\tilde{W}}_{t,T}$-measurable.

We denote
%
%e3.6 #&#
\begin{equation}
\label{aux} X^{t,y}_s=y+B_s-B_t\qquad
\forall t\le s\le T,
\end{equation}
and
consider the following BDSDE:
%
%e3.7 #&#
\begin{eqnarray}
\label{eq0622a} Y^{t,y}_s=F\bigl(X^{t,y}_T
\bigr)+\int^T_s\int_U G
\bigl(a,y,Y^{t,y}_r\bigr)\tilde{W}(\hat{d}r \,da)-\int
^T_sZ^{t,y}_r\,dB_r,
\nonumber
\\[-8pt]
\\[-8pt]
\eqntext{t\le s\le T.}
\end{eqnarray}

BDSDE (\ref{eq0622a}) coincides with BDSDE (\ref{eq0831c}) if we
take $\xi=F(X^{t,y}_T)$ and let the initial time be denoted by $t$
instead of 0 ($t$ is fixed and $s$ varies as shown). We use the
superscript $(t,y)$ to indicate the dependency on the initial state
of the underlying motion.

%th3.4 #&#
\begin{thm}\label{thm0901c}
Suppose that conditions (\ref{eq0831d}) and (\ref{eq0922a})
hold. If the process $\{\tilde{u}_t\}$ is a solution to
(\ref{eq2a6}) such that $\tilde{u}_\cdot\in C([0,T],\cX_1)$ a.s.,
and
%
%e3.8 #&#
\begin{equation}
\label{eq0824b} \EE\int^T_0\|
\tilde{u}_s\|^2_1\,ds<\infty,
\end{equation}
then
\[
\tilde{u}_t(y)=Y^{t,y}_t,
\]
where $Y^{t,y}_s$ is a solution to the BDSDE
(\ref{eq0622a}).
\end{thm}
\begin{pf} Let
%
%e3.9 #&#
\begin{equation}
\label{eq0831b} Y^{t,y}_s=\tilde{u}_s
\bigl(X^{t,y}_s\bigr)\quad\mbox{and}\quad Z^{t,y}_s=
\nabla \tilde{u}_s\bigl(X^{t,y}_s\bigr),\qquad  t\le s
\le T.
\end{equation}
To prove
(\ref{eq0622a}), we need to smooth the function $\tilde{u}_t$. For
any ${\delta}>0$, let
\[
u^{\delta}_t(y)=T_{\delta}\tilde{u}_t(y)\qquad
\forall y\in\RR.
\]
It is well known
that for any $t\ge0$ and ${\delta}>0$, $u^{\delta}_t\in C^\infty$. Applying
$T_{\delta}$ to both sides of (\ref{eq2a6}), we have
%
%e3.10 #&#
\begin{eqnarray}
\label{eq0622c} u^{\delta}_t(y)&=&T_{\delta}F(y)+\int
^T_t\frac12\De u^{\delta}_s(y)\,ds
\nonumber
\\[-8pt]
\\[-8pt]
\nonumber
&&{}+\int^T_t\int_U\int
_\RR p_{\delta}(y-z)G\bigl(a,z,\tilde{u}_s(z)
\bigr)\,dz\tilde{W}(\hat{d}s \,da).
\end{eqnarray}

Let $s=t_0<t_1<\cdots<t_n=T$ be a partition of $[s,T]$. Then
\begin{eqnarray*}
&&u^{\delta}_s\bigl(X^{t,y}_s
\bigr)-T_{\delta}F\bigl(X^{t,y}_T\bigr)
\\
&&\qquad=\sum^{n-1}_{i=0}\bigl(u^{\delta}_{t_i}
\bigl(X^{t,y}_{t_i}\bigr)-u^{\delta
}_{t_i}
\bigl(X^{t,y}_{t_{i+1}}\bigr)\bigr) +\sum
^{n-1}_{i=0}\bigl(u^{\delta}_{t_i}
\bigl(X^{t,y}_{t_{i+1}}\bigr)-u^{\delta
}_{t_{i+1}}
\bigl(X^{t,y}_{t_{i+1}}\bigr)\bigr)
\\
&&\qquad=-\sum^{n-1}_{i=0}\int
^{t_{i+1}}_{t_i}\frac12\De u^{\delta}_{t_i}
\bigl(X^{t,y}_r\bigr)\,dr -\sum^{n-1}_{i=0}
\int^{t_{i+1}}_{t_i}\nabla u^{\delta
}_{t_i}
\bigl(X^{t,y}_r\bigr)\,dB_r
\\
&&\qquad\quad{}+\sum^{n-1}_{i=0}\int
^{t_{i+1}}_{t_i}\frac12\De u^{\delta}_r
\bigl(X^{t,y}_{t_{i+1}}\bigr)\,dr
\\
&&\qquad\quad{}+\sum^{n-1}_{i=0}\int
^{t_{i+1}}_{t_i}\int_\RR\int
_U p_{\delta}\bigl(X^{t,y}_{t_{i+1}}-z
\bigr)G\bigl(a,z,\tilde{u}_r(z)\bigr)\tilde{W}(\hat{d}r \,da)\,dz,
\end{eqnarray*}
where we used It\^o's formula for $u^{\delta}_{t_i}$ (note that
$u^{\delta}_{t_i}$ is independent of $X^{t,y}_r$ and~$B_r$), and SPDE
(\ref{eq0622c}) with $y$ replaced by $X^{t,y}_{t_{i+1}}$. Setting
the mesh size to go to 0, we obtain
%
%e3.11 #&#
\begin{eqnarray}
\label{eq0902e}
&&u^{\delta}_s\bigl(X^{t,y}_s
\bigr)-T_{\delta}F\bigl(X^{t,y}_T\bigr)
\nonumber\\
&&\qquad=-\int^T_s\nabla u^{\delta}_r
\bigl(X^{t,y}_r\bigr)\,dB_r
\\
&&\qquad\quad{} +\int^T_s\int_\RR
\int_U p_{\delta}\bigl(X^{t,y}_r-z
\bigr)G\bigl(a,z,\tilde{u}_r(z)\bigr)\tilde{W}(\hat{d}r \,da)\,dz.
\nonumber
\end{eqnarray}
We take ${\delta}\to0$ on both sides of (\ref{eq0902e}). Note that for
$s>t$,
\begin{eqnarray*}
&&\EE\biggl\llvert \int^T_s\nabla
u^{\delta}_r\bigl(X^{t,y}_r
\bigr)\,dB_r-\int^T_s\nabla
\tilde{u}_r\bigl(X^{t,y}_r\bigr)\,dB_r
\biggr\rrvert^2
\\
&&\qquad=\EE\int^T_s\bigl\llvert \nabla
u^{\delta}_r\bigl(X^{t,y}_r\bigr)-\nabla
\tilde{u}_r\bigl(X^{t,y}_r\bigr)\bigr
\rrvert^2\,dr
\\
&&\qquad\le\EE\int^T_s\int_\RR
\bigl(T_{\delta}\nabla\tilde{u}_r(z)-\nabla
\tilde{u}_r(z)\bigr)^2p_{r-t}(y-z)\,dz \,dr.
\end{eqnarray*}
For $s>t$ fixed, there exists a constant $K_1$, depending on $s-t$,
such that for any $r>s$,
\[
p_{r-t}(y-z)\le K e^{-|y-z|}\le K e^{|y|}e^{-|z|}.
\]
Thus, we may continue the estimate above with
\begin{eqnarray*}
&&\EE\biggl\llvert \int^T_s\nabla
u^{\delta}_r\bigl(X^{t,y}_r
\bigr)\,dB_r-\int^T_s\nabla
\tilde{u}_r\bigl(X^{t,y}_r\bigr)\,dB_r
\biggr\rrvert^2
\\
&&\qquad \le Ke^{|y|}\EE\int^T_s \int
_\RR\bigl(T_{\delta}\nabla\tilde {u}_r(z)-
\nabla \tilde{u}_r(z)\bigr)^2e^{-|z|}\,dz \,dr\\
&&\qquad\to0,
\end{eqnarray*}
where the last step follows from the integrability condition
(\ref{eq0824b}).

The other terms can be estimated similarly. (\ref{eq0622a}) follows
from (\ref{eq0902e}) by taking ${\delta}\to0$.
\end{pf}

%s4 #&#
\section{Uniqueness for SPDE}\label{sec4}

The existence of a solution to SPDE (\ref{eq0623a}) was established
in Section~\ref{sec2}. This section is devoted to the proof of the
uniqueness part of Theorem~\ref{thm0618a}.

\begin{pf*}{Proof of Theorem~\ref{thm0618a} \textup{(Uniqueness)}} Let
$u^j_s, j=1,2,$ be two solutions to SPDE (\ref{eq0623a}). Let $T>0$
be fixed and let $\tilde{u}^j_s=u^j_{T-s}$. Denote\vadjust{\goodbreak}
$u^{j,{\delta}}_s=T_{\delta}\tilde{u}^j_s, j=1, 2$, and let $s>t$
be fixed.
By~(\ref{eq0902e}),
%
%e4.1 #&#
\begin{eqnarray}
\label{eq0902f} &&u^{1,{\delta}}_s\bigl(X^{t,y}_s
\bigr)-u^{2,{\delta}}_s\bigl(X^{t,y}_s\bigr)\nonumber\\
&&\qquad=-\int^T_s\nabla\bigl(u^{1,{\delta}}_s-u^{2,{\delta}}_s
\bigr) \bigl(X^{t,y}_r\bigr)\,dB_r\nonumber
\\[-8pt]
\\[-8pt]
\nonumber
&&\qquad\quad{}+\int^T_s\int_U\int
_\RR p_{\delta}\bigl(X^{t,y}_r-z\bigr) \\
&&\hspace*{89pt}{}\times\bigl(G\bigl(a,z,\tilde{u}^1_r(z)\bigr)-G
\bigl(a,z,\tilde{u}^2_r(z)\bigr)\bigr)\,dz \tilde{W}(
\hat{d}r \,da).\nonumber
\end{eqnarray}
Let $\phi_k$ be
defined as in the proof of Theorem~\ref{thm0901b}. Applying the
It\^o--Pardoux--Peng formula to (\ref{eq0902f}) and $\phi_k$, similarly
to (\ref{eq0902c}), we get
%
%e4.2 #&#
\begin{eqnarray}
\label{eq0904b} &&\EE\phi_k\bigl(u^{1,{\delta}}_s
\bigl(X^{t,y}_s\bigr)-u^{2,{\delta
}}_s
\bigl(X^{t,y}_s\bigr)\bigr)
\nonumber\\
&&\qquad\le\frac12\EE\int^T_s\int
_U\phi''_k
\bigl(u^{1,{\delta
}}_r\bigl(X^{t,y}_r
\bigr)-u^{2,{\delta}}_r\bigl(X^{t,y}_r\bigr)
\bigr)
\nonumber
\\[-8pt]
\\[-8pt]
\nonumber
&&\hspace*{46pt}\qquad\quad{}\times\biggl\llvert \int_\RR p_{\delta}
\bigl(X^{t,y}_r-z\bigr) \bigl(G\bigl(a,z,
\tilde{u}^1_r(z)\bigr)\\
&&\hspace*{106pt}\qquad\quad{}
-G\bigl(a,z,\tilde
{u}^2_r(z)\bigr)\bigr)\,dz\biggr\rrvert^2
\la(da) \,dr.
\nonumber
\end{eqnarray}
Next, we take the limit ${\delta}\to0$ on both sides of (\ref{eq0904b}).
By Lemma~\ref{lem0225a}, $T_{\delta}\tilde{u}^j_s\to\tilde{u}^j_s$ in
$\cX_0$ as ${\delta}\to0$. Taking a subsequence if necessary, we may and
will assume that $T_{\delta}\tilde{u}^j_s(x)\to\tilde{u}^j_s(x)$ for
almost every $x$ with respect to the Lebesgue measure. Therefore,
\[
u^{1,{\delta}}_s\bigl(X^{t,y}_s
\bigr)-u^{2,{\delta}}_s\bigl(X^{t,y}_s\bigr)
\to \tilde{u}^{1}_s\bigl(X^{t,y}_s
\bigr)-\tilde{u}^{2}_s\bigl(X^{t,y}_s
\bigr)\qquad \mbox{a.s.,}
\]
and by the bounded convergence theorem, the left-hand side of
(\ref{eq0904b}) converges to
\[
\EE\phi_k\bigl(\tilde{u}^{1}_s
\bigl(X^{t,y}_s\bigr)-\tilde{u}^{2}_s
\bigl(X^{t,y}_s\bigr)\bigr).
\]

Denote
\[
g_r(a,z)=G\bigl(a,z,\tilde{u}^1_r(z)
\bigr)-G\bigl(a,z,\tilde{u}^2_r(z)\bigr),\qquad (a,z)\in U
\times\RR.
\]
Then, the right-hand side of (\ref{eq0904b}) can be
written as
%
%e4.3 #&#
\begin{eqnarray}
\label{eq0225c} \qquad &&\frac12\EE\int^T_s\int
_\RR\int_U\phi''_k
\bigl(u^{1,{\delta
}}_r(x)-u^{2,{\delta}}_r(x)
\bigr)\bigl|T^U_{\delta} g_r(a,x)\bigr|^2p_{r-t}(x-y)\,dx\,
\la(da)\,dr
\nonumber
\\[-8pt]
\\[-8pt]
\nonumber
&&\qquad=\frac12\EE\int^T_s\bigl\|\bigl(T^U_{\delta}g_r
\bigr) h_r\bigr\|^2_{\cX_0\otimes
L^2(U,\la)}\,dr,
\end{eqnarray}
where $h_r(x)$, $r\ge s$ and
$x\in\RR$, is such that
\[
h_r(x)^2=\phi''_k
\bigl(u^{1,{\delta}}_r(x)-u^{2,{\delta
}}_r(x)
\bigr)e^{|x|}p_{r-t}(x-y).
\]
Note that $h_r(x)$ is bounded by a constant depending on $(k, s-t, y)$. On the other hand,
\[
\|g_r\|^2_{\cX_0\otimes
L^2(U,\la)}\le K\int_\RR
\bigl(1+\bigl|u^1_r(z)\bigr|^2+\bigl|u^2_r(z)\bigr|^2
\bigr)e^{-|z|}\,dz,
\]
which is integrable. By Lemma~\ref{lem0225b} and the dominated
convergence theorem, we see that the limit of the right-hand side of
(\ref{eq0904b}) is equal to
\begin{eqnarray*}
&&\frac12\EE\int^T_s\lim_{{\delta}\to0}
\bigl\|T^U_{\delta}g_r h_r
\bigr\|^2_{\cX_0\otimes L^2(U,\la)}\,dr
\\
&&\qquad=\frac12\EE\int^T_s\|g_r
h_r\|^2_{\cX_0\otimes L^2(U,\la)}\,dr
\\
&&\qquad=\frac12\EE\int^T_s\phi''_k
\bigl(\tilde{u}^{1}_r\bigl(X^{t,y}_r
\bigr)-\tilde {u}^{2}_r\bigl(X^{t,y}_r
\bigr)\bigr) \bigl|\tilde{u}^1_r\bigl(X^{t,y}_r
\bigr)-\tilde{u}^2_r\bigl(X^{t,y}_r
\bigr)\bigr| \,dr.
\end{eqnarray*}

To summarize, we obtain
%
%e4.4 #&#
\begin{eqnarray}
\label{eq0225d} &&\EE\phi_k\bigl(\tilde{u}^{1}_s
\bigl(X^{t,y}_s\bigr)-\tilde{u}^{2}_s
\bigl(X^{t,y}_s\bigr)\bigr)
\nonumber\\
&&\qquad\le\frac12\EE\int^T_s\phi''_k
\bigl(\tilde{u}^{1}_r\bigl(X^{t,y}_r
\bigr)-\tilde {u}^{2}_r\bigl(X^{t,y}_r
\bigr)\bigr) \bigl|\tilde{u}^1_r\bigl(X^{t,y}_r
\bigr)-\tilde{u}^2_r\bigl(X^{t,y}_r
\bigr)\bigr| \,dr
\\
&&\qquad \le k^{-1}T,
\nonumber
\end{eqnarray}
where we used $|z|\phi''_k(z)\le2k^{-1}$ in the last step.

Finally, applying Fatou's lemma for $k\to\infty$, we obtain
\[
\EE\bigl|\tilde{u}^{1}_s\bigl(X^{t,y}_s
\bigr)-\tilde{u}^{2}_s\bigl(X^{t,y}_s
\bigr)\bigr|\le \liminf_{k\to\infty}\EE\phi_k\bigl(\tilde{u}^{1}_s
\bigl(X^{t,y}_s\bigr)-\tilde {u}^{2}_s
\bigl(X^{t,y}_s\bigr)\bigr)\le0.
\]
Therefore, $\tilde{u}^1_s(X^{t,y}_s)-\tilde{u}^2_s(X^{t,y}_s)=0$
a.s. Taking $s\downarrow t$, we get $u^1_t(y)=u^2_t(y)$, a.s.
\end{pf*}

After proving the pathwise (strong) uniqueness and weak existence of
the solution for SPDE (\ref{eq0623a}), we verify its (weak)
uniqueness. For finite dimensional It\^o equations, Yamada and
Watanabe~\cite{YW} proved that weak existence and strong uniqueness
imply strong existence and weak uniqueness. Kurtz~\cite{Kur}
considered this problem in an abstract setting. To apply Kurtz's
result to SPDE (\ref{eq0623a}), we convert it to an SPDE driven by a
sequence of independent Brownian motions. Let $\{h_j\}^\infty_{j=1}$
be a CONS of $L^2(U,\cU,\la)$, and define
\[
B^j_t=\int^t_0\int
_U h_j(a)W(ds \,da),\qquad j=1,2,\ldots.
\]
Letting $B_t=(B^j_t)^\infty_{j=1}$, it is easy to see that
(\ref{eq0623a}) is equivalent to the following SPDE:
%
%e4.5 #&#
\begin{equation}
\label{eq0623a'} {u}_t(y)=F(y)+\sum
^\infty_{j=1}\int^t_0
G_j\bigl(y,{u}_s(y)\bigr)\,dB^j_s+
\int^t_0\frac12\De{u}_s(y)\,ds,
\end{equation}
where
\[
G_j(y,u)=\int_U G(a,y,u)h_j(a)
\la(da).
\]

Denote
\[
S_1=\CC\bigl([0,T],\cX_0\bigr)\quad\mbox{and}\quad S_2=\CC\bigl([0,T],\RR^\infty\bigr).
\]
Let $\{f_k\}^\infty_{k=1}\subset C^\infty_0(\RR)$ be a dense subset
of $\cX_0$ and $\Ga\dvtx  S_1\times S_2\to\RR$ be the measurable
functional defined by
\[
\Ga(u_\cdot,B_\cdot)=\sum^\infty_{k=1}
\sup_{t\le T}\bigl|{\gamma }^k_t\bigr|\wedge
2^{-k},
\]
where
\[
\Ga^k_t=\langle{u}_t,f_k
\rangle-\langle F,f_k\rangle-\int^t_0
\biggl\langle {u}_s,\frac12\De f_k\biggr\rangle \,ds-\sum
^\infty_{j=1}\int^t_0
\int_\RR G_j\bigl(y,{u}_s(y)
\bigr)f(y)\,dy \,dB^j_s.
\]
Then, SPDE (\ref{eq0623a'}) can be rewritten as
\[
\Ga(u_\cdot,B_\cdot)=0.
\]

The following theorem is a direct consequence of Proposition 2.10 in
Kurtz~\cite{Kur}, which is needed for next section.

%th4.1 #&#
\begin{thm}\label{thm0309a}
If $(u^i_\cdot)$, $i=1,2,$ are two solutions of SPDE (\ref{eq0623a})
(may be defined on different stochastic bases) such that
\[
\EE\sup_{t\le T}\bigl\|u^i_t\bigr\|^2_0<
\infty,\qquad  i=1,2,
\]
then their laws in $\CC([0,T],\cX_0)$ coincide.
\end{thm}

%s5 #&#
\section{Measure-valued processes}\label{sec5}

In this section, we give the proofs of three applications of
Theorem~\ref{thm0618a} to measure-valued processes.

Recall
that SBM $\mu_t$ is defined as the unique solution to the following
martingale problem (MP): $\forall f\in C^2_b(\RR)$, the process
%
%e5.1 #&#
\begin{equation}
\label{eq0902a} M^f_t\equiv\mu_t(f)-\mu(f)-
\int^t_0\mu_s\biggl(\frac12
f''\biggr)\,ds
\end{equation}
is a continuous square-integrable martingale
with
%
%e5.2 #&#
\begin{equation}
\label{eq0902b} \bigl\langle M^f\bigr\rangle_t=\int
^t_0\mu_s\bigl(f^2
\bigr)\,ds.
\end{equation}

Now, we present:

\begin{pf*}{Proof of Theorem~\ref{thm4sbm}} Suppose that $\mu_t$ is an SBM
and $u_t$ is defined by~(\ref{eq0831a}). Let $f\in C^2_0(\RR)$ and
$g(y)=\int^\infty_yf(x)\,dx$. Then
%
%e5.3 #&#
\begin{eqnarray}
\label{eq0825b} \langle u_t,f\rangle&=&\mu_t(g)
\nonumber\\
&=&\mu_0(g)+\int^t_0
\mu_s\biggl(\frac12 g''
\biggr)\,ds+M^g_t
\\
&=&\langle F,f\rangle+\int^t_0\biggl\langle
u_s,\frac12 f''\biggr\rangle
\,ds+M^g_t.
\nonumber
\end{eqnarray}
Let $\mathcal{S}'(\RR)$ be the space of Schwartz distributions and
define the $\mathcal{S}'(\RR)$-valued process $N_t$ by
$N_t(f)=M^g_t$ for any $f\in C^\infty_0(\RR)$. Then, $N_t$ is an
$\mathcal{S}'(\RR)$-valued continuous square-integrable martingale
with
\begin{eqnarray*}
\bigl\langle N(f)\bigr\rangle_t&=&\int^t_0
\int_\RR g(y)^2\mu_s(dy)\,ds
\\
&=&\int^t_0\int_\RR g
\bigl(u_s^{-1}(a)\bigr)^2\,da \,ds
\\
&=&\int^t_0\int_\RR
\biggl(\int_\RR1_{a\le u_s(y)}f(y)\,dy\biggr)^2\,da
\,ds,
\end{eqnarray*}
where $u_s^{-1}$ is the generalized inverse of the nondecreasing
function $u_s$, that is,
\[
u_s^{-1}(a)=\sup \bigl\{x\in\RR\dvtx  u_s(x)<a
\bigr\}.
\]

Let ${\gamma}\dvtx \RR_+\times\Om\to L_{(2)}(H,H)$ be defined as
\[
{\gamma}(s,{\omega})f(a)=\int_\RR1_{a\le u_s(x)}f(x)\,dx\qquad
\forall f\in H,
\]
where $H=L^2(\RR)$ and $L_{(2)}(H,H)$ is the space consisting of all
Hilbert--Schmidt operators on $H$. By Theorem 3.3.5 of Kallianpur and
Xiong~\cite{KX1}, on an extension of the original stochastic basis,
there exists an $H$-cylindric Brownian motion $B_t$ such that
\[
N_t(f)=\int^t_0 \bigl\langle{
\gamma}(s,\omega)f,\,dB_s\bigr\rangle_H.
\]
Let $\{h_j\}$ be a CONS of the Hilbert space $H$ and define random
measure $W$ on $\RR_+\times\RR$ as
\[
W\bigl([0,t]\times A\bigr)=\sum^\infty_{j=1}
\langle1_A,h_j\rangle B^{h_j}_t.
\]
It is easy to show that $W$ is a Gaussian white noise random measure
on $\RR_+\times\RR$ with intensity $ds \,da$. Furthermore,
\[
N_t(f)=\int^t_0\int
_\RR\int_\RR1_{a\le u_s(x)}f(x)\,dx
W(ds \,da).
\]
Plugging back to (\ref{eq0825b}) verifies that $u_t$ is a solution
to (\ref{eq0609a}).

On the other hand, suppose that $\{u_t\}$ is a weak solution to SPDE
(\ref{eq0609a}) with $F\in\cX_0$ being nondecreasing. Let $\mu_0$
be the measure determined by $F$. Let $\nu_t$ be an SBM with initial
$\mu_0$. Define the function-valued process $\hat{u}_t$ by
\[
\hat{u}_t(y)=\int^y_0
\nu_t(dx)\qquad \forall y\in\RR.
\]
By the above result, $\hat{u}_t$ is a solution to SPDE (\ref{eq0609a})
with initial $F$. Here we remark that (\ref{eq0609a}) coincides with
(\ref{eq0623a}) if we take $U=\RR$, $\la(da)=da$ and
$G(a,y,u)=1_{0\le a\le u}+1_{u\le a\le0}$. By the weak uniqueness
(Theorem~\ref{thm0309a}) of the solution to this SPDE, $(u_t)$ and
$(\hat{u}_t)$ have the same distribution, implying $(\mu_t)$ and
$(\nu_t)$ have the same distribution. This proves that $(\mu_t)$ is
an SBM.
\end{pf*}

The result for the Fleming--Viot process is similar so we only
provide a sketch.

\begin{pf*}{Sketch of the proof of Theorem~\ref{thm4fv}} The uniqueness of
SPDE (\ref{eq0609b}) follows from Theorem~\ref{thm0618a} by taking
$U=[0,1]$, $\la(da)=da$ and $G(a,y,u)=1_{0\le a\le u}-u$.

Suppose that
$\{u_t\}$ is a weak solution to the SPDE (\ref{eq0609b}), and
$\{\mu_t\}$ is defined by~(\ref{eq0831a}). Then for any $f\in
C^3_0(\RR)$,
\begin{eqnarray*}
\mu_t(f)&=&-\bigl\langle u_t,f'\bigr
\rangle
\\
&=&-\bigl\langle F,f'\bigr\rangle-\int^t_0
\int_\RR\frac12 u_s(y)f'''(y)\,dy
\,ds
\\
&&{}-\int_\RR\int^t_0\int
^1_0\bigl(1_{a\le
u_s(y)}-u_s(y)
\bigr)W(ds \,da)f'(y)\,dy
\\
&=&\mu(f)+\int^t_0\mu_s\biggl(
\frac12f''\biggr)\,ds
\\
&&{}+\int^t_0\int^1_0
\bigl(f\bigl(u_s^{-1}(a)\bigr)-\mu_s(f)
\bigr)W(ds \,da).
\end{eqnarray*}
Thus
\begin{eqnarray*}
N^f_t&\equiv&\mu_t(f)-\mu(f)-\int
^t_0\mu_s\biggl(
\frac12f''\biggr)\,ds
\\
&=&\int^t_0\int^1_0
\bigl(f\bigl(u_s^{-1}(a)\bigr)-\mu_s(f)
\bigr)W(ds \,da)
\end{eqnarray*}
is a continuous square-integrable martingale with
\begin{eqnarray*}
\bigl\langle N^f\bigr\rangle_t&=&\int
^t_0\int^1_0
\bigl(f\bigl(u_s^{-1}(a)\bigr)-\mu_s(f)
\bigr)^2\,da \,ds
\\
&=&\int^t_0\bigl(\mu_s
\bigl(f^2\bigr)-\mu_s(f)^2\bigr)\,ds.
\end{eqnarray*}
The proof of other direction is similar, so we omit it.
\end{pf*}

Finally, we present:

\begin{pf*}{Proof of Theorem~\ref{thm0904a}} Denote by $\HH$ the
reproducing kernel Hilbert space (RKHS) of the covariance function
$\phi$. In other words, $\HH$ is the completion of the linear span
of the functions $\{\phi(x,\cdot)\dvtx x\in\RR\}$ with respect to the
inner product
\[
\bigl\langle\phi(x,\cdot),\phi(y,\cdot)\bigr\rangle_\HH=\phi(x,y).
\]
We refer the reader to Kallianpur~\cite{Kal}, page 139, for more
details on RKHS. Let $\{h_j\}$ be a CONS of $\HH$. Let $U=\NN$ and
let $\la(da)$ be the counting measure. Note that
\begin{eqnarray*}
\phi(x,y)&=&\sum^\infty_{j=1}\bigl\langle
\phi(x,\cdot),h_j\bigr\rangle_\HH \bigl\langle\phi(y,
\cdot),h_j\bigr\rangle_\HH
\\
&=&\int_U\rho(a,x)\rho(a,y)\la(da),
\end{eqnarray*}
where $\rho(a,x)=\langle\phi(x,\cdot),h_a\rangle_\HH$.

Let $\cS(\RR)$ be the space of rapidly decreasing functions on $\RR
$; cf. Definition~1.3.4 in Kallianpur and Xiong~\cite{KX} for its
definition. For any $h\in\cS(\RR)$, we define
\[
B_t(h)=\int^t_0\int
_\RR h(x)B(x,ds)\,dx.
\]
Then, $B_t$ is an $\cS'(\RR)$-valued martingale with
\begin{eqnarray*}
\bigl\langle B(h)\bigr\rangle_t&=&\int^t_0
\int_\RR\int_\RR h(x)h(y)\phi (x,y)\,dx
\,dy \,ds
\\
&=&\sum^\infty_{j=1}\int^t_0
\biggl\llvert \int_\RR h(x)\rho(j,x)\,dx\biggr
\rrvert^2\,ds.
\end{eqnarray*}
Analogously to the proof of Theorem~\ref{thm4sbm}, there exists a
sequence of independent Brownian motions $B^j_t$ such that
\[
B_t(h)=\sum^\infty_{j=1}\int
^t_0 \int_\RR h(x)
\rho(j,x)\,dx \,dB^j_s.
\]
Let
\[
W\bigl([0,t]\times\{j\}\bigr)=W^j_t,\qquad  j=1,2,\ldots.
\]
Then, $W$ is a space--time white noise random measure on $\RR_+\times
U$ with intensity $dt\, \la(da)$, and
\[
B(x,dt)=\int_U\rho(a,x)W(dt \,da).
\]
Let
\[
G(a,y,u)=\rho(a,y)\sqrt{u}.
\]
Then, (\ref{eq0901a}) is a special case of SPDE (\ref{eq0623a}) and
conditions (\ref{eq0831d}) and (\ref{eq0922a}) are satisfied. The
conclusion of Theorem~\ref{thm0904a} then follows from Theorem
\ref{thm0618a}.
\end{pf*}

\section*{Acknowledgments}
I would like to thank a number of individuals for their help on this
article. Zenghu Li and Shige Peng provided helpful conversations
about this research. Don Dawson, Carl Mueller, Ed Perkins, Xinghua
Zheng and Xiaowen Zhou gave useful comments after a preliminary
version was circulated. I also gratefully acknowledge two anonymous
referees for many constructive suggestions which improved this paper
substantially.

% imsref loaded by akundreckaite, 2013-01-07 12:21:04
% imsref loaded by akundreckaite, 2013-01-07 12:38:46
%
% imsref loaded by akundreckaite, 2013-01-07 15:06:24

%suskaldyti doi

\printaddresses

\end{document}